\documentclass[12pt]{article}

\usepackage{cite}
\usepackage{amsfonts, amssymb, amsmath, amsthm, epsf}

\numberwithin{figure}{section}

\theoremstyle{plain}
\newtheorem{theorem}{Theorem}

\theoremstyle{definition}
\newtheorem{remark}{Remark}

\newcommand{\Dom}{{\rm D}}

\newcommand{\const}{{\rm const}}

\renewcommand{\phi}{{\varphi}}

\newcommand{\cK}{{\mathcal K}}

\newcommand{\cO}{{\mathcal O}}

\newcommand{\cR}{{\mathcal R}}

\newcommand{\pG}{{\partial G}}

\newcommand{\oG}{{\overline G}}

\newcommand{\bP}{{\mathbf P}}

\newcommand{\bI}{{\mathbf I}}

\newcommand{\bbR}{{\mathbb R}}

\begin{document}

\begin{center}
{\bf\large On nonexistence of Feller semigroups in the nontransversal case}

\smallskip

Pavel Gurevich\footnote{This research was supported by Russian
Foundation for Basic Research (project No.~07-01-00268) and the
Alexander von Humboldt Foundation.}
\end{center}

\begin{abstract}
We give three examples of second-order elliptic operators with
nonlocal boundary conditions of the Ventsel type that admit a
closure in the space of continuous functions, but do not generate
a Feller semigroup (i.e., a strongly continuous contractive
nonnegative semigroup).
\end{abstract}

The question of existence of a strongly continuous contractive nonnegative semigroup
(Feller semigroup) of operators acting between the spaces of continuous functions arises
in the theory of Markov processes. Feller semigroups describe (from the probabilistic
viewpoint) the motion of a Markov particle in a region. A general form of a generator of
such a semigroup on an interval was investigated in~\cite{Feller2}.  In the
multidimensional case, it was proved that the generator of a Feller semigroup is an
elliptic differential operator (possibly with degeneration) whose domain of definition
consists of continuous functions satisfying nonlocal conditions which involve an integral
over the closure of the region with respect to a nonnegative Borel
measure~\cite{Ventsel}. The inverse problem remains open: given an elliptic
integro-differential operator whose domain of definition is described by nonlocal
boundary conditions, whether or not the closure of this operator is a generator of a
Feller semigroup.   The order of nonlocal terms is less than the order of local terms in
the {\it transversal\/} case~\cite{SU, BCP, Taira3,  Ishikawa, GalSkJDE}, and these
orders coincide in the more difficult {\it nontransversal\/} case~\cite{GalSkJDE} (see
also the bibliography therein).

In~\cite{SkRJMP95}, it was given an example of nonlocal operator (involving a
transformation of the boundary into itself)  whose closure is not a generator of a Feller
semigroup. Here we give three examples of nonexistence of Feller semigroups in the cases
where transformations $\Omega(y)$ (in {\it nontransversal} nonlocal conditions) map the
boundary inside the region. For any $y$ on the boundary, the above Borel measure is the
delta function supported at the point $\Omega(y)$ from the closure of the region. We note
that Conditions~3.3 and 3.6 from~\cite{GalSkJDE} fail in our first and second examples
and Conditions 3.5 and 3.9 from~\cite{GalSkJDE} fail in our third example.

{\bf 1. ``Jumps'' with nonzero probability to outside of a neighborhood of
process-termination points.} Let $G\subset\bbR^2$ be a bounded region with smooth
boundary $\pG=\Gamma_1\cup\Gamma_2\cup\cK$, where $\Gamma_1$ and $\Gamma_2$ are
$C^\infty$ curves open and connected in the topology of $\pG$ such that
$\Gamma_1\cap\Gamma_2=\varnothing$ and $\overline{\Gamma_1}\cap\overline{\Gamma_2}=\cK$;
the set $\cK$ consists of two points $g_1$ and $g_2$. We assume that the region $G$
coincides with the plane angle of opening $\pi$ in some $\varepsilon$-neighborhood
$\cO_\varepsilon(g_i)$ of the point $g_i$, $i=1,2$.

Consider the following nonlocal condition:
\begin{equation}
u(y)-b_1(y)u(\Omega_1(y)) =0,\  y\in\Gamma_1;\qquad u(y) =0, \
y\in\overline{\Gamma_2},\label{eqEx1_2-3}
\end{equation}
where $b_1\in C^\infty(\overline{\Gamma_1})$, $0\le b_1(y)\le 1$,
$b_1(y)=\const>0$ for $y\in\cO_{\varepsilon/2}(g_1)$, $b_1(y)=0$
for $y\notin\cO_{ \varepsilon}(g_1)$;  $\Omega_1$ is a smooth
nondegenerate transformation defined on a neighborhood
 of the curve $\overline{\Gamma_1}$,
$\Omega_1(\Gamma_1)\subset G$, $\Omega_1(g_1)\in G$, and $\Omega_1(y)$ is the composition
of  rotation about the point $g_1$ and shift by some vector for
$y\in\cO_\varepsilon(g_1)$. Probabilistically, the Dirichlet condition means that the
Markov particle  is absorbed (the process terminates) as it reaches the boundary at the
point $y\in\overline{\Gamma_2}$; the nonlocal condition  means that the particle
``jumps'' from the point $y\in\Gamma_1$ to the point $\Omega_1(y)\in G$ with probability
$b_1(y)$ after some random time.

We consider the unbounded operator $\bP_1: \Dom(\bP_1)\subset
C_1(\overline G)\to C_1(\overline G)$ given by
\begin{equation*}
\bP_1 u=\Delta u,\qquad u\in \Dom(\bP_1)=\{u\in C_1(\oG): \Delta
u\in C_1(\overline G)\},
\end{equation*}
where $C_1(\overline G)$ is the set of functions from $C(\oG)$
satisfying nonlocal conditions~\eqref{eqEx1_2-3}, $\Delta$ is the
Laplace operator acting in the sense of distributions.

{\bf 2. ``Jumps'' out from the conjugation points that are not
process-termination points.}   Consider the following nonlocal
condition:
\begin{equation}
u(y)-b_1(y)u(\Omega_1(y)) =0,\ y\in\overline{\Gamma_1};\qquad
u(y)-b_2(y)u(\Omega_2(y)) =0, \  y\in\Gamma_2,\label{eqEx3_2-3}
\end{equation}
where $b_j\in C^\infty(\overline{\Gamma_j})$, $0\le b_j(y)\le 1$,
$b_j(y)=b^*>0$ for $y\in\cO_{\varepsilon/2}(g_1)$, $b_j(y)=0$ for
$y\notin\cO_{\varepsilon}(g_1)$;  $\Omega_j$  is a smooth
nondegenerate transformation defined on a neighborhood  of the
curve $\overline{\Gamma_j}$, $\Omega_j(\Gamma_j)\subset G$,
$\Omega_j(g_1)\in G$,    $\Omega_1(g_1)\ne \Omega_2(g_1)$, and
$\Omega_j(y)$ is the composition of   rotation about the point $g_1$ and
shift by some vector for $y\in\cO_\varepsilon(g_1)$.

We consider the unbounded operator $\bP_2: \Dom(\bP_2)\subset
C_2(\overline G)\to C_2(\overline G)$ given by
\begin{equation*}
\bP_2 u=\Delta u,\qquad u\in \Dom(\bP_2)=\{u\in C_2(\oG): \Delta
u\in C_2(\overline G)\},
\end{equation*}
where $C_2(\overline G)$ is the set of functions from $C(\oG)$
satisfying nonlocal conditions~\eqref{eqEx3_2-3}.

{\bf 3.  ``Jumps'' with probability one within  a neighborhood of the
process-termination points.}   Consider the
following nonlocal condition:
\begin{equation}
u(y)-b_j(y)u(\Omega_j(y))=0,\  y\in\Gamma_j,\ j=1,2;\qquad u(y)
=0,\ y\in\cK,\label{eqEx2_2-3}
\end{equation}
where  $b_j\in C^\infty(\overline{\Gamma_j})$, $0\le b_j(y)\le 1$,
$b_j(y)=1$ for $y\in\cO_{\varepsilon/2}(g_1)$, $b_j(y)=0$ for
$y\notin\cO_{\varepsilon}(g_1)$;   $\Omega_j$ is a smooth
nondegenerate transformation defined on a neighborhood of the
curve $\overline{\Gamma_j}$, $\Omega_j(\Gamma_j)\subset G$,
$\Omega_j(g_1)=g_1$, and $\Omega_j(y)$ is the rotation by the angle
$\pi/2$ inwards the region $G$ for $y\in\cO_\varepsilon(g_1)$.

We consider the unbounded operator $\bP_3: \Dom(\bP_3)\subset
C_3(\overline G)\to C_3(\overline G)$ given by
\begin{equation*}
\bP_3 u=\Delta u,\qquad u\in \Dom(\bP_3)=\{u\in C_3(\oG): \Delta
u\in C_3(\overline G)\},
\end{equation*}
where $C_3(\overline G)$ is the set of functions from $C(\oG)$
satisfying nonlocal conditions~\eqref{eqEx2_2-3}.

\begin{theorem}\label{thEx1-3}
The operators ${\bP_j}$  admit the closure
$\overline{\bP_j}:\Dom(\overline{\bP_j})\subset C_j(\oG)\to
C_j(\oG)$ {\rm (}$j=1,2,3${\rm )}. The operators
$\overline{\bP_j}$ {\rm (}$j=1,2,3${\rm )} are  not generators of
Feller semigroups.
\end{theorem}

\begin{remark}
One can prove that $C_j(\oG)\setminus\overline{\cR({\bP_j}-q\bf I)}\ne\varnothing$ for
sufficiently small $q>0$. Hence,
$C_j(\oG)\setminus{\cR(\overline{\bP}_j-q\bI)}\ne\varnothing$. Combining this fact with
the Hille--Iosida theorem, we obtain Theorem~\ref{thEx1-3}.
\end{remark}

The author is grateful to Prof. A.L. Skubachevskii for attention to this work.

\end{document}